\documentclass[11pt]{article}
\usepackage{latexsym,amssymb}

\def\demo{\noindent{\bf Proof. }}
\def\QED{\hfill$\Box$}
\newtheorem{Theorem}{Theorem}[section]
\newtheorem{Lemma}[Theorem]{Lemma}
\newtheorem{Corollary}[Theorem]{Corollary}
\newtheorem{Proposition}[Theorem]{Proposition}
\newtheorem{Remark}[Theorem]{Remark}
\newtheorem{Example}[Theorem]{Example}
\newtheorem{Conjecture}[Theorem]{Conjecture}
\newtheorem{Definition}[Theorem]{Definition}

\begin{document}
\topmargin3mm
\hoffset=-1cm
\voffset=-1.5cm
\

\medskip

\begin{center}
{\large\bf A note on Rees algebras and the MFMC property }
\vspace{6mm}\\
%\footnotetext{2000 {\it Mathematics Subject 
%Classification}. 13F55, 90C10, 13F20, 90C27, 52B20.} 
%\footnotetext{{\it Key words\/}.  
%Rees cone, max-flow min-cut property, Hilbert bases, 
%covering polyhedron.}

\medskip

Isidoro Gitler\footnote{Partially supported by CONACyT grants
49835-F, 49251-F and 
SNI.},\,  Carlos E. Valencia 
\\
and
\\
\smallskip
Rafael H. Villarreal
\\
\medskip
{\small Departamento de Matem\'aticas}\vspace{-1mm}\\ 
{\small Centro de Investigaci\'on y de Estudios Avanzados del
IPN}\vspace{-1mm}\\   
{\small Apartado Postal 14--740}\vspace{-1mm}\\ 
{\small 07000 M\'exico City, D.F.}\vspace{-1mm}\\ 
{\small e-mail: {\tt vila@math.cinvestav.mx}}\vspace{3mm}
\end{center}
\date{}

\begin{abstract} 
\noindent We study irreducible representations of Rees cones and 
characterize the max-flow min-cut property 
of clutters in terms of the normality of Rees algebras 
and the integrality of certain polyhedra. Then we present some
applications to combinatorial optimization and commutative algebra. 
As a byproduct we obtain an effective method, based on
the program {\it Normaliz\/} \cite{B}, to determine whether a given
clutter satisfy the max-flow min-cut property. Let $\cal C$ be a clutter 
and let $I$ be its edge ideal. We prove 
that $\cal C$ has the max-flow min-cut property if and only if 
$I$ is normally torsion free, that is, $I^i=I^{(i)}$ 
for all $i\geq 1$, where $I^{(i)}$ is the $i${\it th} symbolic power
of $I$.  
\end{abstract}

\section{Introduction}\label{intro}

Let $R=K[x_1,\ldots,x_n]$ be a polynomial ring over a field $K$ and
let $I\subset R$ be a monomial ideal
minimally generated by $x^{v_1},\ldots,x^{v_q}$. As usual we will use 
$x^a$ as an abbreviation for $x_1^{a_1} \cdots x_n^{a_n}$, 
where $a=(a_1,\ldots,a_n)\in \mathbb{N}^n$. Consider the 
$n\times q$ matrix $A$ with column vectors $v_1,\ldots,v_q$. A {\it
clutter\/} with 
vertex set  
$X$ is a family of subsets of $X$, called edges, none 
of which is included in another. A basic example  
of clutter is a graph. If 
$A$ has entries in $\{0,1\}$, then $A$ defines in a 
natural way a {\it clutter\/} ${\cal C}$ by taking
$X=\{x_1,\ldots,x_n\}$ as vertex set and  
$E=\{S_1,\ldots,S_q\}$ as edge set, where
$S_i$ is the support of $x^{v_i}$, i.e., the set of variables
that occur in $x^{v_i}$. In this case we call $I$ 
the {\it edge ideal\/} of the 
clutter ${\cal C}$ and write $I=I({\cal C})$. Edge ideals
are also called {\it facet ideals} \cite{faridi}. This notion has
been studied by Faridi \cite{faridi1} and Zheng \cite{zheng}. The
matrix $A$ is often 
refer to as the {\it incidence\/} matrix of $\cal C$.

The {\it Rees algebra} of $I$ 
is the $R$-subalgebra:
$$
R[It]:=R[\{x^{v_1}t,\ldots,x^{v_q}t\}]\subset R[t],
$$
where $t$ is a new variable. In our situation $R[It]$ is also a
$K$-subalgebra of $K[x_1,\ldots,x_n,t]$. The {\it Rees cone\/} of $I$
is the rational polyhedral  
cone in $\mathbb{R}^{n+1}$, denoted by $\mathbb{R}_+{\cal A}'$, consisting 
of the non-negative linear combinations of the set 
$${\cal
A}':=\{e_1,\ldots,e_n,(v_1,1),\ldots,(v_q,1)\}\subset\mathbb{R}^{n+1},$$  
where $e_i$ is the $i${\it th} unit vector. Thus ${\cal A}'$
is the set of exponent vectors of the set of monomials 
$\{x_1,\ldots,x_n,x^{v_1}t,\ldots,x^{v_q}t\}$, that generate $R[It]$ 
as a $K$-algebra. 

The first main result of this note (Theorem~\ref{irre}) shows that 
the irreducible 
representation of the Rees cone, as a finite intersection of closed
half-spaces, can be  
expressed essentially in terms of the vertices of the {\it 
set covering polyhedron\/}: 
\[
Q(A):=\{x \in \mathbb{R}^n  \, | \; x \geq{0}, \;
xA\geq{\mathbf 1} \}.
\]
Here ${\mathbf 1}=(1,\ldots,1)$. The second main result
(Theorem~\ref{mfmc=i+n}) is an 
algebro-combinatorial 
description of the max-flow min-cut property of the clutter $\cal C$
in terms of a purely  
algebraic property (the normality of $R[It]$) and an integer programming 
property (the integrality of the rational polyhedron $Q(A)$). Some 
applications will be shown. For instance we give an effective method,
based on the 
program {\it Normaliz\/} \cite{B}, to  determine whether a given
clutter satisfy the max-flow min-cut property
(Remark~\ref{mfmc-effe-method}). We prove that $\cal C$ has
the max-flow min-cut property if and only if $I^i=I^{(i)}$ 
for $i\geq 1$, where $I^{(i)}$ is the $i${\it th} symbolic power
of $I$ (Corollary~\ref{mfmc-characterization}). There are other
interesting links between algebraic
properties of Rees algebras and combinatorial optimization
problems of 
clutters \cite{reesclu}. 

Our main references for Rees algebras and combinatorial optimization
are \cite{BSV,Vas} and \cite{korte} respectively.

\section{Preliminaries}\label{prelim}

For convenience we quickly recall some basic results, terminology,
and notation from  polyhedral geometry. 

A set 
$C\subset\mathbb{R}^n$ is a {\it polyhedral set\/} (resp. {\it cone})
if 
$C=\{x\vert\, Bx\leq b\}$ for some matrix $B$ and some vector $b$
(resp. $b=0$).  By the finite basis theorem 
\cite[Theorem 4.1.1]{webster} a polyhedral cone $C\subsetneq$ $\mathbb{R}^n$ 
has two representations:
\begin{description}
\item[\it Minkowski representation]  $C=\mathbb{R}_+{\cal B}$ with 
${\cal B}=\{\beta_1,\ldots,\beta_r\}$ a finite set, and
\item[\it Implicit representation] $C=H_{c_1}^+\cap\cdots\cap H_{c_s}^+$ for 
some $c_1,\ldots,c_s\in\mathbb{R}^n\setminus\{0\}$,  
\end{description}
where $\mathbb{R}_+$ is the set
of non-negative real numbers, $\mathbb{R}_+{\cal B}$ is the cone 
generated by $\cal B$ consisting of the set of linear combinations of 
$\cal B$ with coefficients in $\mathbb{R}_+$, 
$H_{c_i}$ is the hyperplane of ${\mathbb R}^n$ through the origin  
with normal vector ${c_i}$, and $H_{c_i}^+=\{x\vert\, 
\langle x, c_i\rangle\geq 0\}$ 
is the positive {\em closed half-space\/} bounded by $H_{c_i}$. Here 
$\langle\, ,\, \rangle$ denotes the usual inner product.  These 
two representations satisfy the {\it duality theorem for cones\/}: 
\begin{equation}\label{duality}
H_{\beta_1}^+\cap\cdots\cap H_{\beta_r}^+=\mathbb{R}_+c_1+\cdots+\mathbb{R}_+c_s,
\end{equation}
see \cite[Corollary~7.1a]{Schr} and its proof.
The {\it dual cone} of 
$C$ is defined as
$$
C^*:=\bigcap_{c\in C} H_c^+=\bigcap_{a\in{\cal B}}H_a^+.
$$
By the duality theorem $C^{**}=C$. An implicit 
representation of $C$ is called {\it irreducible\/} if none of the 
closed half-spaces $H_{c_1}^+,\ldots,H_{c_s}^+$ can be omitted from
the intersection.  
Note that the left hand side of Eq.~(\ref{duality}) is an 
irreducible representation of $C^*$ if and only if no proper subset 
of ${\cal B}$ generates $C$.

\section{Rees cones, normality and the MFMC property}\label{irres}

To avoid repetions, throughout the rest of this note we keep the
notation and assumptions  
of Section~\ref{intro}. 

Notice that the Rees cone $\mathbb{R}_+{\cal A}'$ has dimension
$n+1$. 
A subset $F\subset \mathbb{R}^{n+1}$ is called a 
{\em facet\/} of $\mathbb{R}_{+}{\cal A}'$ if  $F=\mathbb{R}_+{\cal
A}'\cap H_a$ for some hyperplane $H_a$ such that $\mathbb{R}_+{\cal
A}'\subset H_a^{+}$ and $\dim(F)=n$.  
It is not hard to see 
that the set
$$
F=\mathbb{R}_+{\cal A}'\cap H_{e_i}\ \ \ \ (1\leq i\leq n+1)
$$
defines a facet of $\mathbb{R}_+{\cal A}'$ if and only if either $i=n+1$ 
or $1\leq i\leq n$ and $\langle e_i,v_j\rangle=0$ for some column
$v_j$ of $A$. Consider  
the index set
$${\cal J}=\{1\leq i\leq n | \; \langle e_i,v_j\rangle=0 \mbox{ for
some } j\}\cup \{n+1\}. 
$$
Using \cite[Theorem~3.2.1]{webster} it is seen that the Rees cone has
a unique irreducible  
representation 
\begin{equation}\label{okayama-car} 
{\mathbb R}_+{\cal A}'=\left(\bigcap_{i\in \cal J}H_{e_i}^+\right)\bigcap\left( 
\bigcap_{i=1}^r H_{a_i}^+\right)
\end{equation}
such that $0\neq a_i\in\mathbb{Q}^{n+1}$ and 
$\langle a_i,e_{n+1}\rangle=-1$ for all $i$. A point $x_0$ 
is called a {\it vertex\/} or an {\it extreme point\/} of $Q(A)$ if
$\{x_0\}$ is a proper face of $Q(A)$. 

\begin{Lemma}\label{jul13-03} Let $a=(a_{i1},\ldots,a_{iq})$ be the
$i${\it th} row of  
the matrix $A$ and define $k=\min\{a_{ij}\vert\, 1\leq j\leq q\}$. If 
$a_{ij}>0$ for all $j$, then $e_i/k$ is a vertex of $Q(A)$. 
\end{Lemma}

\demo Set $x_0=e_i/k$. Clearly $x_0\in Q(A)$ and $\langle x_0,v_j\rangle=1$ for 
some $j$. Since $\langle x_0,e_\ell\rangle=0$ for $\ell\neq i$, the point 
$x_0$ is a basic feasible solution of $Q(A)$. Then by
\cite[Theorem~2.3]{bertsimas} $x_0$ is a  
vertex of $Q(A)$. \QED 

\begin{Theorem}\label{irre} Let $V$ be the 
vertex set of $Q(A)$. Then 
$$
\mathbb{R}_+{\cal A}'=\left(\bigcap_{i\in {\cal J}} 
H^+_{e_i}\right)\bigcap \left(\bigcap_{\alpha\in V}H^+_{(\alpha,-1)}\right)
$$
is the irreducible 
representation of the Rees cone of $I$.
\end{Theorem}

\demo Let $V=\{\alpha_1,\ldots,\alpha_p\}$ be the set of vertices of $Q(A)$ and let 
$${\cal B}=\{e_i\vert\, i\in{\cal J}\}\cup\{(\alpha,-1)\vert\, \alpha\in V\}.$$
First we dualize Eq.~(\ref{okayama-car}) and use the duality theorem for cones to obtain
\begin{eqnarray}
(\mathbb{R}_+{\cal A}')^\ast&=&\{y\in\mathbb{R}^{n+1}\vert\, \langle y,
x\rangle\geq 0,\, \forall\, x\in\mathbb{R}_+{\cal A}'\}\nonumber\\  
&=&H_{e_1}^+\cap\cdots\cap 
H_{e_{n}}^+\cap H_{(v_1,1)}^+\cap\cdots\cap H_{(v_q,1)}^+\nonumber\\
 &=& \sum_{i\in \cal J}\mathbb{R}_+e_i+
\mathbb{R}_+a_1+\cdots+\mathbb{R}_+a_r.\label{jul5-03}
\end{eqnarray}
Next we show the equality
\begin{equation}\label{valenciaeq}
(\mathbb{R}_+{\cal A}')^\ast=\mathbb{R}_+{\cal B}.
\end{equation}
The right hand side is clearly contained in the left hand side
because a vector  
$\alpha$ belongs to $Q(A)$ if and only if $(\alpha,-1)$ is in
$(\mathbb{R}_+{\cal A}')^\ast$.  
To prove the reverse containment observe that by Eq.~(\ref{jul5-03})
it suffices to show  
that $a_k\in\mathbb{R}_+{\cal B}$ for all $k$. Writing $a_k=(c_k,-1)$ and using 
$a_k\in(\mathbb{R}_+{\cal A}')^\ast$ gives $c_k\in Q(A)$. The set
covering polyhedron  
can be written as
$$
Q(A)=\mathbb{R}_+e_1+\cdots+\mathbb{R}_+e_n+{\rm conv}(V),
$$
where ${\rm conv}(V)$ denotes the convex hull of $V,$ this follows
from the structure  
of polyhedra by noticing 
that the characteristic cone of $Q(A)$ is precisely $\mathbb{R}_+^n$ (see 
\cite[Chapter~8]{Schr}). 
Thus we can write 
$$
c_k=\lambda_1e_1+\cdots+\lambda_ne_n+\mu_1\alpha_1+\cdots+\mu_p\alpha_p,
$$
where $\lambda_i\geq 0$, $\mu_j\geq 0$ for all $i,j$ and $\mu_1+\cdots+\mu_p=1$. If 
$1\leq i\leq n$ and $i\notin{\cal J}$, then the $i${\it th} row of $A$ has all its 
entries positive. Thus by Lemma~\ref{jul13-03} we get that $e_i/k_i$ is a 
vertex of $Q(A)$ for some 
$k_i>0$.  
To avoid cumbersome notation we denote $e_i$ and $(e_i,0)$ simply by $e_i$, from 
the context the meaning of $e_i$ should be clear. Therefore from the equalities
$$
\sum_{i\notin \cal J}\lambda_ie_i=\sum_{i\notin \cal
J}\lambda_ik_i\left(\frac{e_i}{k_i}\right) 
=\sum_{i\notin \cal
J}\lambda_ik_i\left(\frac{e_i}{k_i},-1\right)+\left(\sum_{i\notin \cal J} 
\lambda_ik_i\right)e_{n+1}
$$
we conclude that $\sum_{i\notin \cal J}\lambda_ie_i$ is in
$\mathbb{R}_+{\cal B}$.  
From the identities
\begin{eqnarray*}
a_k&=&(c_k,-1)=\lambda_1e_1+\cdots+\lambda_ne_n+\mu_1(\alpha_1,-1)+
\cdots+\mu_p(\alpha_p,-1)\\
&=&\sum_{i\notin \cal J}\lambda_ie_i+\sum_{i\in{\cal J}\setminus\{n+1\}}
\hspace{-4mm}\lambda_ie_i+
\sum_{i=1}^p\mu_i(\alpha_i,-1)
\end{eqnarray*}
we obtain that $a_k\in\mathbb{R}_+{\cal B}$, as required. Taking 
duals in Eq.~(\ref{valenciaeq}) 
we get
\begin{equation}\label{valenciaeq1}
\mathbb{R}_+{\cal A}'=\bigcap_{a\in\cal B} H_a^+.
\end{equation}
Thus, by the comments at the end of Section~\ref{prelim}, the proof 
reduces to showing that 
$\beta\notin\mathbb{R}_+({\cal B}\setminus\{\beta\})$ 
for all $\beta\in{\cal B}$. To prove this 
we will assume that $\beta\in\mathbb{R}_+({\cal B}\setminus\{\beta\})$ for 
some $\beta\in{\cal B}$ 
and derive a contradiction. 

Case (I): $\beta=(\alpha_{j},-1)$. For simplicity assume
$\beta=(\alpha_p,-1)$. We can  
write 
$$(\alpha_p,-1)=\sum_{i\in \cal J}\lambda_ie_i+
\sum_{j=1}^{p-1}\mu_j(\alpha_j,-1),\ \ \ \ \ 
(\lambda_i\geq0;\mu_j\geq 0).
$$
Consequently
\begin{eqnarray}
\alpha_p&=&\sum_{i\in {\cal J}\setminus\{n+1\}}\hspace{-4mm}\lambda_ie_i+
\sum_{j=1}^{p-1}\mu_j\alpha_j\label{jul14-03}\\ 
-1&=&\lambda_{n+1}-(\mu_1+\cdots+\mu_{p-1}).\label{jul14-1-03}
\end{eqnarray}
To derive a contradiction we claim that 
$Q(A)=\mathbb{R}_+^n+{\rm conv}(\alpha_1,\ldots,\alpha_{p-1})$, which 
is impossible because by \cite[Theorem~7.2]{Bron} the vertices of
$Q(A)$ would be  
contained in $\{\alpha_1,\ldots,\alpha_{p-1}\}$. To prove the claim
note that the right hand  
side is clearly contained in the left hand side. For the other
inclusion take 
$\gamma\in Q(A)$ and write 
\begin{eqnarray*}
\gamma&=&\sum_{i=1}^nb_ie_i+\sum_{i=1}^{p}c_i\alpha_i\ \ \ \ \ \ \  
(b_i,c_i\geq 0;\sum_{i=1}^pc_i=1)\\ 
&\stackrel{(\ref{jul14-03})}{=}
&\delta+\sum_{i=1}^{p-1}(c_i+c_p\mu_i)\alpha_i\ \ \ \ \ \ \ \
(\delta\in\mathbb{R}_+^n). 
\end{eqnarray*}
Therefore using the inequality
$$
\sum_{i=1}^{p-1}(c_i+c_p\mu_i)=\sum_{i=1}^{p-1}c_i+c_p
\left(\sum_{i=1}^{p-1}\mu_i\right)\stackrel{(\ref{jul14-1-03})}{=}
(1-c_p)+c_p(1+\lambda_{n+1})\geq 1
$$
we get $\gamma\in\mathbb{R}_+^n+{\rm
conv}(\alpha_1,\ldots,\alpha_{p-1})$. This proves  
the claim.

Case (II): $\beta=e_k$ for some $k\in{\cal J}$. First we consider the subcase 
$k\leq n$. The subcase $k=n+1$ can be treated similarly. We can write
$$
e_k=\sum_{i\in {\cal J}\setminus\{k\}}\hspace{-4mm}\lambda_ie_i+
\sum_{i=1}^{p}\mu_i(\alpha_i,-1),\ \ \ \ \ 
(\lambda_i\geq0;\mu_i\geq 0).
$$ 
From this equality we get $e_k=\sum_{i=1}^p\mu_i\alpha_i$. 
Hence $e_kA\geq(\sum_{i=1}^p\mu_i)\mathbf{1}>0$, a contradiction 
because $k\in\cal J$ and $\langle e_k,v_j\rangle=0$ 
for some $j$. \QED

\paragraph{Clutters with the max-flow min-cut property} For the rest of
this section we assume that $A$ is a $\{0,1\}$-matrix, i.e., $I$ is a
square-free monomial ideal. 

\begin{Definition}\rm The clutter
$\cal C$ has the {\it max-flow min-cut\/} (MFMC) 
property if both sides 
of the LP-duality equation
\begin{equation}\label{jun6-2-03}
{\rm min}\{\langle \alpha,x\rangle \vert\, x\geq 0; xA\geq \mathbf{1}\}=
{\rm max}\{\langle y,\mathbf{1}\rangle \vert\, y\geq 0; Ay\leq\alpha\} 
\end{equation}
have integral optimum solutions $x$ and $y$ for each non-negative
integral vector $\alpha$.  
\end{Definition}

It follows from \cite[pp.~311-312]{Schr} that $\cal C$ has the MFMC
property if and only  
if the maximum in Eq.~(\ref{jun6-2-03}) has an optimal integral
solution $y$ for each  
non-negative integral vector $\alpha$. In optimization terms
\cite{korte} this means that the clutter $\cal C$ has the MFMC
property if and only if 
the system of linear inequalities $x\geq 0;\ xA\geq\mathbf{1}$ that
define $Q(A)$ is {\it totally dual integral\/} (TDI).  The polyhedron
$Q(A)$ is said 
to be {\it integral\/} if $Q(A)$ has only integral vertices. 

Next we recall two descriptions of the integral closure of 
$R[It]$ that yield some formulations of the normality
property of $R[It]$. 
Let
$\mathbb{N}{\cal A}'$ be the subsemigroup of  
$\mathbb{N}^{n+1}$ generated by ${\cal A}'$, consisting of the 
linear combinations of ${\cal A}'$ with non-negative integer
coefficients. The Rees algebra of the ideal $I$ can be written as
\begin{eqnarray}
R[It]&=&
K[\{x^at^b\vert\, (a,b)\in\mathbb{N}{\cal A}'\}]\label{may6-06-2}\\ &=&
R\oplus It\oplus\cdots\oplus I^{i}t^i\oplus\cdots
\subset R[t].\label{may6-06-3}
\end{eqnarray}
According to \cite[Theorem~7.2.28]{monalg} and \cite[p.~168]{Vas1} the  
integral closure of $R[It]$ 
in its field of fractions can be expressed as
\begin{eqnarray}
\overline{R[It]}&=&K[\{x^at^b\vert\, (a,b)\in
\mathbb{Z}{\cal A}'\cap \mathbb{R}_+{\cal A}'\}]\label{may6-06}\\ 
&=&R\oplus
\overline{I}t\oplus\cdots\oplus
\overline{I^i}t^i\oplus\cdots,\label{jun05-1-03} 
\end{eqnarray}
where $\overline{I^i}=(\{x^a\in R\vert\, \exists\, p\geq
1;(x^a)^{p}\in I^{pi}\})$ is the integral 
closure of $I^i$ and $\mathbb{Z}{\cal A}'$ is the subgroup of $\mathbb{Z}^{n+1}$ 
generated by ${\cal A}'$. Notice that in our situation we have the
equality $\mathbb{Z}{\cal A}'=\mathbb{Z}^{n+1}$. Hence, by
Eqs.~(\ref{may6-06-2}) to 
$(\ref{jun05-1-03})$, we get that $R[It]$ is a normal domain if and
only if any of the following two conditions hold: 
(a) $\mathbb{N}{\cal A}'=
\mathbb{Z}^{n+1}\cap\mathbb{R}_+{\cal A}'$, 
(b) $I^{i}=\overline{I^i}$ for $i\geq 1$.

\begin{Theorem}\label{mfmc=i+n} The clutter $\cal C$ has 
the {\rm MFMC} property if and only if $Q(A)$ is an integral
polyhedron and $R[It]$ is a normal domain. 
\end{Theorem}
\demo $\Rightarrow)$ By \cite[Corollary~22.1c]{Schr} the polyhedron $Q(A)$ is 
integral. Next we show that $R[It]$ is normal. Take
$x^\alpha t^{\alpha_{n+1}}\in\overline{R[It]}$. Then 
$(\alpha,\alpha_{n+1})\in {\mathbb Z}^{n+1}\cap \mathbb{R}_+{\cal
A}'$. Hence  
$Ay\leq {\alpha}$ and $\langle y,\mathbf{1}\rangle=\alpha_{n+1}$ for some vector 
$y\geq 0$. Therefore one concludes that the optimal value of the linear program 
$$
\max\{\langle y,\mathbf{1}\rangle\vert\; y \geq \mathbf{0};\; Ay\leq{\alpha}\}
$$
is greater or equal than $\alpha_{n+1}$. Since $A$ has the MFMC
property, this linear  
program has an optimal integral solution $y_0$. Thus there exists an
integral vector $y_0'$ 
such that 
$$
\mathbf{0} \leq y_0' \leq y_0\ \mbox{ and }\ |y_0'|=\alpha_{n+1}.
$$
Therefore
$$
\left(\hspace{-2mm}\begin{array}{c}
\alpha\\ \alpha_{n+1}
\end{array}\hspace{-1mm}
\right)=\left(\hspace{-2mm}\begin{array}{c}
A\\ \mathbf{1}
\end{array}\hspace{-1mm}
\right)y_0'+
\left(\hspace{-2mm}\begin{array}{c}
A\\ \mathbf{0}
\end{array}\hspace{-1mm}
\right)(y_0-y_0')+\left(\hspace{-2mm}\begin{array}{c}
\alpha\\ 0
\end{array}\hspace{-1mm}
\right)-\left(\hspace{-2mm}\begin{array}{c}
A\\ \mathbf{0}
\end{array}\hspace{-1mm}
\right)y_0
$$
and $(\alpha,\alpha_{n+1})\in \mathbb{N}{\cal A}'$. This proves that
$x^\alpha t^{\alpha_{n+1}}\in{R[It]}$, as required.

$\Leftarrow)$ Assume that $A$ does not satisfy the MFMC 
property. There exists an $\alpha_0\in {\mathbb N}^n$ such that if ${y}_0$ 
is an optimal solution of the linear program:
$$
\max\{\langle y,\mathbf{1}\rangle\vert\;y\geq \mathbf{0};\;
Ay\leq{\alpha}_0\},\eqno(*) 
$$
then ${y}_0$ is not integral. We claim that also the optimal value 
$|y_0|=\langle{y}_0,\mathbf{1}\rangle$ of this linear program is not integral. If 
$|y_0|$ is integral, 
then $(\alpha_0,|y_0|)$ is in $\mathbb{Z}^{n+1}\cap\mathbb{R}_+{\cal
A}'$. As 
$R[It]$ is normal, we get that $(\alpha_0,|y_0|)$ is in 
$\mathbb{N}{\cal A}'$, but this readily yields that the linear program 
$(*)$ has an integral optimal solution, a contradiction. This completes the 
proof of the claim. 

Now, consider the dual linear program:
$$
\min\{\langle x,\alpha_0\rangle\vert \; x \geq \mathbf{0},\;  xA\geq \mathbf{1}\}.
$$
By \cite[Theorem~4.1.6]{webster}) the optimal value of this linear
program is attained at a  
vertex $x_0$ of $Q(A)$. Then by the LP duality
theorem~\cite[Theorem~3.16 ]{korte} we 
get $\langle x_0,\alpha_0\rangle=|{y}_0|\notin {\mathbb Z}$. Hence $x_0$ is not 
integral, a contradiction to the integrality of the 
set covering polyhedron $Q(A)$. \QED

\begin{Remark}\label{mfmc-effe-method}\rm 
The program {\it Normaliz\/} \cite{B,BK} computes 
the irreducible representation of a Rees cone and the integral
closure of $R[It]$. Thus one can effectively use 
Theorems~\ref{irre} and \ref{mfmc=i+n} to determine whether a given
clutter $\cal C$ as the max-flow min-cut property. See example below
for a simple illustration.
\end{Remark}

\begin{Example}\label{sept20-3-03}\rm Let 
$I=(x_1x_5,x_2x_4,x_3x_4x_5,x_1x_2x_3)$. Using {\it
Normaliz\/} \cite{B} with the input file: 
\begin{small}
\begin{verbatim}
4
5
1 0 0 0 1
0 1 0 1 0
0 0 1 1 1
1 1 1 0 0
3
\end{verbatim}
\end{small}
we get the output file:
\begin{small}
\begin{verbatim}
9 generators of integral closure of Rees algebra: 
  1  0  0  0  0  0
  0  1  0  0  0  0
  0  0  1  0  0  0
  0  0  0  1  0  0
  0  0  0  0  1  0
  1  0  0  0  1  1
  0  1  0  1  0  1
  0  0  1  1  1  1
  1  1  1  0  0  1

10 support hyperplanes: 
   0   0   1   1   1  -1
   1   0   0   0   0   0
   0   1   0   0   0   0
   0   0   0   0   0   1
   0   0   1   0   0   0
   1   0   0   1   0  -1
   0   0   0   1   0   0
   0   0   0   0   1   0
   0   1   0   0   1  -1
   1   1   1   0   0  -1
\end{verbatim}
\end{small}
\noindent The first block shows the exponent vectors of the generators
of the integral
closure of $R[It]$, thus $R[It]$ is normal.  The second block 
shows the irreducible representation of the Rees cone of $I$, thus
using Theorem~\ref{irre} we obtain that $Q(A)$ is integral. 
Altogether Theorem~\ref{mfmc=i+n} proves that the 
clutter $\cal C$ associated to $I$ has the max-flow min-cut property.
\end{Example}

\begin{Definition}\rm A set $C\subset X$ is a 
{\it minimal vertex cover\/} of a clutter $\cal C$ if every edge of $\cal C$ 
contains at least one vertex in $C$ and $C$ is minimal w.r.t. 
this property. A set of edges of $\cal C$ is {\it
independent} if no 
two of them have a common vertex.  We denote by ${\alpha}_0({\cal
C})$ the smallest number of vertices in any  
minimal vertex cover of $\cal C$, and by $\beta_1({\cal C})$ the 
maximum number of independent edges of ${\cal C}$.  
\end{Definition}

\begin{Definition}\rm Let $X=\{x_1,\ldots,x_n\}$ and let 
$X'=\{x_{i_1},\ldots,x_{i_r},x_{j_1},\ldots,x_{j_s}\}$ be a subset of
$X$. A {\it minor\/} of $I$ is a proper ideal $I'$ of
$R'=K[X\setminus X']$ obtained from  
$I$ by making 
$x_{i_k}=0$ and $x_{j_\ell}=1$ for all $k,\ell$. The ideal $I$ is
considered itself a minor.  
A {\it minor\/} of $\cal C$ 
is a clutter ${\cal C}'$ that corresponds to a minor $I'$.  
\end{Definition}

Recall that a ring is called 
{\it reduced} if $0$ is its only nilpotent element. The {\it associated graded
ring\/} of $I$ is the quotient ring
${\rm gr}_I(R):=R[It]/IR[It]$.

\begin{Corollary}\label{konig} If the associated graded ring ${\rm gr}_I(R)$ is 
reduced, then 
$\alpha_0({\cal C}')=\beta_1({\cal C}')$ for any minor ${\cal C}'$ of
$\cal C$. 
\end{Corollary}

\demo As the reducedness of ${\rm gr}_I(R)$ is preserved if we 
make a variable $x_i$ equal to $0$ or $1$, we may assume that 
${\cal C}'={\cal C}$. From \cite[Proposition~3.4]{normali} and
Theorem~\ref{irre} it follows  
that the ring  ${\rm gr}_I(R)$ is reduced if and only if $R[It]$ 
is normal and $Q(A)$ is integral. Hence by Theorem~\ref{mfmc=i+n} we obtain that 
the LP-duality equation
$$
{\rm min}\{\langle\mathbf{1},x\rangle \vert\, x\geq 0; xA\geq \mathbf{1}\}=
{\rm max}\{\langle y,\mathbf{1}\rangle \vert\, y\geq 0; Ay\leq\mathbf{1}\}
$$
has optimum integral solutions $x$, $y$. To complete the 
proof notice that the left hand side of this 
equality is $\alpha_0({\cal C})$ and the right hand side 
is $\beta_1({\cal C})$. \QED

\medskip

Next we state an algebraic version of a conjecture 
\cite[Conjecture~1.6]{cornu-book} which to our best knowledge 
is still open:  

\begin{Conjecture}\label{conforti-cornuejols}\rm 
If $\alpha_0({\cal C}')=\beta_1({\cal C}')$ for all 
minors $\cal C'$ of $\cal C$, then the associated graded ring ${\rm
gr}_I(R)$ 
is reduced.
\end{Conjecture} 

\begin{Proposition} Let $B$ be the matrix with column 
vectors $(v_1,1),\ldots,(v_q,1)$. If $x^{v_1},\ldots,x^{v_q}$ are monomials of the
same degree $d\geq 2$ and ${\rm gr}_I(R)$ is reduced, then $B$
diagonalizes over $\mathbb{Z}$ to an identity matrix.
\end{Proposition}

\demo As $R[It]$ is normal, the result follows from 
\cite[Theorem~3.9]{ehrhart}. \QED

\medskip

This result suggest the following weaker conjecture:  

\begin{Conjecture}[\rm Villarreal]\label{git-val-vi}\rm Let $A$ be a
$\{0,1\}$-matrix  
such that the number of $1$'s in every column of $A$ has a constant value 
$d\geq 2$. If $\alpha_0({\cal C}')=\beta_1({\cal C}')$ for all 
minors $\cal C'$ of $\cal C$, then the quotient group 
$\mathbb{Z}^{n+1}/((v_1,1),\ldots,(v_q,1))$ is torsion-free.
\end{Conjecture} 

\paragraph{Symbolic Rees algebras} Let
$\mathfrak{p}_1,\ldots,\mathfrak{p}_s$ be the minimal primes  
of the edge ideal $I=I({\cal C})$ and let 
$C_k=\{x_i\vert\, x_i\in\mathfrak{p}_k\}$, for $k=1,\ldots,s$, 
be the corresponding minimal vertex covers of the clutter $\cal C$. We set 
$$
\ell_k=(\textstyle\sum_{x_i\in C_k}e_i,-1)  \ \ \ \ (k=1,\ldots,s).
$$

The {\it symbolic Rees algebra\/} of $I$ is the $K$-subalgebra:
$$
R_s(I)=R+I^{(1)}t+I^{(2)}t^2+\cdots+I^{(i)}t^i+\cdots\subset R[t],
$$
where $I^{(i)}=\mathfrak{p}_1^i\cap\cdots\cap\mathfrak{p}_s^i$ is the
$i${\it th\/} symbolic power of $I$. 

\begin{Corollary}\label{oct30-03} The following conditions are
equivalent
\begin{description}
\item{\rm (a)} $Q(A)$ is integral.\vspace{-1mm}
\item{\rm (b)} ${\mathbb R}_+{\cal A}'=H_{e_1}^+\cap \cdots\cap H_{e_{n+1}}^+
\cap H_{\ell_1}^+\cap\cdots\cap H_{\ell_s}^+$.
\vspace{-1mm}
\item{\rm (c)} $\overline{R[It]}=R_s(I)$, i.e., 
$\overline{I^i}=I^{(i)}$ for all $i\geq 1$.\vspace{-1mm}
\end{description}
\end{Corollary}

\demo The integral vertices of $Q(A)$ are precisely the 
vectors $a_1,\ldots,a_s$, where 
$a_k=\sum_{x_i\in C_k}e_i$ for $k=1,\ldots,s$. Hence by
Theorem~\ref{irre} we obtain that (a) is equivalent to (b). 
By \cite[Corollary~3.8]{normali} we get that (b) is 
equivalent to (c). \QED

\begin{Corollary}\label{mfmc-characterization} Let $\cal C$ be a
clutter and 
let $I$ be its edge
ideal. Then $\cal C$ has the max-flow min-cut property if and
only if $I^i=I^{(i)}$ for all $i\geq 1$.  
\end{Corollary}

\demo It follows at once from Corollary~\ref{oct30-03} and
Theorem~\ref{mfmc=i+n}. \QED


\begin{thebibliography}{50}


\bibitem{bertsimas} D. Bertsimas and J. N. Tsitsiklis, {\it Introduction 
to linear optimization\/}, Athena Scientific, Massachusetts, 1997.

\bibitem{Bron}{A. Br\o ndsted, {\it Introduction to Convex
Polytopes},  Graduate
Texts in  Mathematics {\bf 90}, Springer-Verlag, 1983.}

\bibitem{BSV}{P. Brumatti, A. Simis and W.~V. Vasconcelos, Normal
{R}ees algebras, J.~Algebra {\bf 112} (1988), 26--48.}

\bibitem{B}{W. Bruns and R. Koch, {\it Normaliz} -- a program 
for computing normalizations of affine semigroups, 1998. 
Available via anonymous ftp from 
{\rm ftp.mathematik.Uni-Osnabrueck.DE/pub/osm/kommalg/software}.}

\bibitem{BK}{W. Bruns and R. Koch, Computing the integral 
closure of an affine semigroup. Effective methods in 
algebraic and analytic geometry, 2000 (Krak\'ow). Univ. 
Iagel. Acta Math. {\bf 39} (2001), 59--70.}


\bibitem{cornu-book}{G. Cornu\'ejols, {\it Combinatorial optimization: 
Packing and covering\/}, CBMS-NSF Regional Conference Series in Applied 
Mathematics {\bf 74}, SIAM (2001).}

\bibitem{ehrhart}{C. Escobar, J.
Mart\'\i nez-Bernal and R. H. Villarreal, Relative volumes and minors
in monomial subrings, Linear Algebra Appl. {\bf 374} (2003), 275--290.}

\bibitem{normali} C. Escobar, R. H. Villarreal and Y. Yoshino, Torsion
freeness and normality of blowup rings of monomial ideals, 
{\it Commutative Algebra\/}, Lect. Notes Pure Appl. Math. 
{\bf 244}, Chapman \& Hall/CRC, Boca Raton, FL, 2006, pp. 69-84. 

\bibitem {faridi} S. Faridi, The facet ideal of a simplicial complex, 
Manuscripta Math. {\bf 109} (2002), 159-174.

\bibitem{faridi1} S. Faridi, Cohen-Macaulay properties of square-free 
monomial ideals, J. Combin. Theory Ser. A {\bf 109(2)} (2005),
299-329. 


\bibitem{reesclu}{I. Gitler, E. Reyes and R. H. Villarreal, 
Blowup algebras of square--free monomial ideals and some links to
combinatorial optimization problems, Rocky Mountain J. Math., to appear.} 

\bibitem{korte} B. Korte and J. Vygen, 
{\it Combinatorial Optimization Theory 
and Algorithms\/}, Springer-Verlag, 2000.

\bibitem{Schr}{A. Schrijver, {\it Theory of Linear and 
Integer Programming\/}, John Wiley \& Sons, New York, 1986.}

\bibitem{Vas}{W. V. Vasconcelos, {\it Arithmetic of Blowup 
Algebras\/}, London Math. Soc., Lecture Note Series 
{\bf 195}, Cambridge
University Press, Cambridge, 1994. }

\bibitem{Vas1}{W. V. Vasconcelos, {\it Computational Methods in
Commutative Algebra and Algebraic Geometry\/}, 
Springer-Verlag, 1998.}

\bibitem{monalg}{R. H. Villarreal, {\it Monomial Algebras\/},
Monographs and Textbooks  
in Pure and Applied Mathematics {\bf 238}, Marcel Dekker, New York, 2001.}

\bibitem{webster}{R. Webster, {\it Convexity\/}, Oxford University
Press, Oxford, 1994.} 

\bibitem{zheng} X. Zheng, Resolutions of facet ideals, 
Comm. Algebra {\bf 32(6)} (2004),  2301--2324. 

\end{thebibliography}
\end{document}